\documentclass{article}
\usepackage{amssymb, amsmath}
\usepackage[matrix, arrow, curve]{xy}
\newcommand \lap {\lambda^{\prime}}
\newcommand \pip {\pi^{\prime}}

\newcommand \la {\lambda}
\newcommand \gz {{\cal G}}

\newcommand \R {{\cal R}}
\newcommand \A {{\cal A}}

\newcommand \q {{\bf q}}
\newcommand \Prob {\mathbb P}

\newcommand \U {{\cal U}}

\newcommand \Y {{\cal Y}}

\newcommand \F {{\cal F}}

\newcommand \wab {{W}_{{\cal A},B}}

\newcommand \omabz {{\Omega}_{{\cal A},B}^{{\mathbb Z}}}
\newcommand \omab { {\Omega}_{{\cal A},B}}

\newtheorem {lemma} {Lemma}
\newtheorem {theorem} {Theorem}
\newtheorem {proposition} {Proposition}

\begin{document}

\title{Exponential decay of correlations for the Rauzy-Veech-Zorich 
induction map.}
\author{Artur Avila\footnote{CNRS UMR 7599, Laboratoire de Probabilit\'es et
Mod\`eles al\'eatoires,
Universit\'e Pierre et Marie Curie--Bo\^\i te courrier 188,
75252--Paris Cedex 05, France} \ and 
Alexander Bufetov\footnote{Department of Mathematics, Rice University}} 
\maketitle

\begin{abstract}

We prove exponential mixing
for the Rauzy-Veech-Zorich induction map on the space of
interval exchange transformations (Theorem \ref{main}).

\end{abstract}

\section{Introduction.}

This section mainly serves to fix notation.
We follow the notation and the conventions of \cite{bufetov}.

\subsection{Rauzy operations $a$ and $b$.}

Let $\pi$ be a permutation on $m$ symbols. The permutation $\pi$ will 
always be assumed
irreducible in the sense that $\pi\{1,\dots ,k\}=\{1,\dots,k\}$ iff $k=m$.
Rauzy operations $a$ and $b$ are defined by the formulas:

$$
a\pi(j)=\begin{cases}
\pi j,&\text{if $j\leq \pi^{-1}m$;}\\
                 \pi m,&\text{if $j=\pi^{-1}m+1$;}\\
                  \pi(j-1),&\text{ other $j$.}
\end{cases}
$$

$$
b\pi(j)=\begin{cases}
\pi j,&\text{if $\pi j\leq \pi m$;}\\
                 \pi j+1,&\text{if $\pi m<\pi j<m$;}\\
                  \pi m+1,&\text{ if $\pi j=m$.}
\end{cases}
$$

These operations preserve irreducibility.
The {\it Rauzy class} of a permutation
$\pi$ is defined as the set of all permutations that can be obtained
from $\pi$ by repeated application of the operations $a$ and $b$.

For $i,j=1,\dots,m$,
denote by $E_{ij}$ an $m\times m$ matrix of which
the $i,j$-th element is equal to $1$, all others to $0$.
Let $E$ be the $m\times m$-identity matrix.

Following Veech \cite{veech}, we introduce the matrices
$$
A(a, \pi)=\sum_{i=1}^{\pi^{-1}(m)}E_{ii}+E_{m, \pi^{-1}m+1}+
\sum_{i=\pi^{-1}m+1}^m E_{i,i+1},
$$
$$
A(b, \pi)=E+E_{m,\pi^{-1}m}
$$

For a vector $\la\in{\mathbb R}^m$, $\la=(\la_1, \dots, \la_m)$,
we write
$$
|\la|=\sum_i \la_i.
$$
Let $\Delta_{m-1}$ be the unit simplex in ${\mathbb R}^m$:
$$
\Delta_{m-1}=\{\la\in {\mathbb R}_+^m: |\la|=1\}.
$$
{\it The space} $\Delta(\R)$ {\it
of interval exchange transformations}, corresponding
to a Rauzy class $\R$, is defined by the formula
$$
\Delta(\R)=\Delta_{m-1}\times \R.
$$
Denote
$$
\Delta_{\pi}^+=\{\la\in\Delta_{m-1}| \ \la_{\pi^{-1}m}>\la_m\},\ \ 
\Delta_{\pi}^-=\{\la\in\Delta_{m-1}| \ \la_m>\la_{\pi^{-1}m}\},
$$
and
$$
\Delta^+ = \cup_{\pip\in{\cal R}(\pi)}\Delta_{\pip}^+,\ \ 
\Delta^- = \cup_{\pip\in{\cal R}(\pi)}\Delta_{\pip}^-.
$$

{\it The Rauzy-Veech induction} is a map
$$
{\cal T}: \Delta(\R)\to \Delta(\R),
$$
defined by the formula
$$
{\cal T}(\la,\pi)=\begin{cases}
(\frac{A(\pi, a)^{-1}\la}{|A(\pi, a)^{-1}\la|}, a\pi), &\text{if $\la\in\Delta^-$;}\\
(\frac{A(\pi, b)^{-1}\la}{|A(\pi, b)^{-1}\la|}, b\pi), &\text{if $\la\in\Delta^+$.}
\end{cases}
$$
Veech \cite{veech} showed that the Rauzy-Veech induction
has an absolutely continuous ergodic invariant measure on $\Delta(\R)$; that measure is,
however, infinite.

Following Zorich \cite{zorich}, define the function $n(\la,\pi)$ in the following way:
$$
n(\la, \pi)=\begin{cases}
\inf \{k>0:{\cal T}^k(\la,\pi)\in\Delta^-\},&\text{if $\la\in\Delta_{\pi}^+$;}\\
  \inf \{k>0: {\cal T}^k(\la,\pi)\in\Delta^+\},&\text{if $\la\in\Delta_{\pi}^-$.}
\end{cases}
$$

{The Rauzy-Veech-Zorich induction} is defined by the formula
$$
{\cal G}(\la,\pi)={\cal T}^{n(\la,\pi)}(\la,\pi).
$$

\begin{theorem}[Zorich\cite{zorich}]
\label{zorichthm}
The map ${\cal G}$ has an ergodic invariant probability measure, absolutely
continuous with respect to the Lebesgue measure class on $\Delta(\R)$.
\end{theorem}

This invariant measure will be denoted by $\nu$.

\subsection{Symbolic dynamics for the induction map.}

This subsection briefly describes the symbolic dynamics for the map $\gz$.
The notation follows \cite{bufetov}.

Consider the alphabet
$$
{\cal A}=\{ (c,n,\pi) | \  c=a \ {\rm or}\  b, n\in {\mathbb N}\}.
$$

For $w_1\in\A$, $w_1=(c_1, n_1, \pi_1)$, we write $c_1=c(w_1), \pi_1=\pi(w_1), n_1=n(w_1)$.
For $w_1,w_2\in\A$, $w_1=(c_1, n_1, \pi_1)$, $w_2=(c_2, n_2, \pi_2)$, define
the function $B(w_1, w_2)$ in the following way:
$B(w_1, w_2)=1$ if $c_1^{n_1}\pi_1=\pi_2$ and $c_1\neq c_2$ and $B(w_1, w_2)=0$ otherwise.

Introduce the space of words
$$
{\cal W}_{{\cal A}, B}=\{w=w_1\dots w_n | \
w_i\in{\cal A}, B(w_i, w_{i+1})=1 \  {\rm for \ all}\  i=1, \dots,n\}.
$$

For a word $w\in\wab$, we denote by $|w|$ its length, i.e., the number
of symbols in it;
given two words $w(1), w(2)\in\wab$, we denote by $w(1)w(2)$ their
concatenation.
Note that the word $w(1)w(2)$  need not
belong to $\wab$, unless a compatibility condition is satisfied by
the last symbol of $w(1)$ and the first symbol of $w(2)$.

To each word assign the corresponding renormalization matrix as follows. 
For $w_1\in\A$, $w_1=(c_1, n_1, \pi_1)$, set 
$$
A(w_1)=A(c_1, \pi_1)A(c_1, c_1\pi_1)\dots A(c_1, c_1^{n_1-1}\pi_1),
$$
and for $w\in W_{\A,B}$, $w=w_1\dots w_n$, set 
$$
A(w)=A(w_1)\dots A(w_n).
$$

Words from $\wab$ act on permutations from $\R$: namely, 
if $w_1\in \A$, $w_1=(c_1, n_1, \pi_1)$, then we set $w_1\pi_1=c_1^{n_1}\pi_1$. 
For permutations $\pi\neq \pi_1$, the symbol $w_1\pi$ is not defined. 
Furthermore, for 
$w\in W_{\A,B}$, $w=w_1\dots w_n$, we write
$$
w\pi=w_n(w_{n-1}(\dots w_1\pi)\dots ),
$$
assuming the right-hand side of the expression is defined. 
Finally, if $\pip=w\pi$, then we also write $\pi=w^{-1}\pip$. 

We say that $w_1\in\A$ is compatible with $(\la,\pi)\in\Delta(\R)$ if 

\begin{enumerate}
\item either $\la\in\Delta_{\pi}^{+}$, $c_1=a$, and $a^{n_1}\pi_1=\pi$
\item or $\la\in\Delta_{\pi}^{-}$, $c_1=b$, and $b^{n_1}\pi_1=\pi$.
\end{enumerate}

We say that a word $w\in W_{\A,B}$, $w=w_1\dots w_n$ is compatible with $(\la,\pi)$ if 
$w_n$ is compatible with $(\la,\pi)$.  We shall also sometimes say that $(\la,\pi)$ is 
compatible with $w$ instead of saying that $w$ is compatible with $(\la,\pi)$.


Now consider the sequence spaces 
$$
\Omega_{\A, B}=\{\omega=\omega_1\dots \omega_n\dots |\  \omega_n\in\A, \ B(\omega_n, \omega_{n+1})=1  \ {\rm for \ all}\  n\in {\mathbb N}\},
$$
and
$$
\Omega_{\A, B}^{\mathbb Z}=\{\omega=\dots \omega_{-n}\dots \omega_1\dots \omega_n\dots |\  \omega_n\in\A, \ 
B(\omega_n, \omega_{n+1})=1
\ {\rm for \ all}\  n\in{\mathbb Z}\}.
$$

Denote by $\sigma$ the right shift on both these spaces.
By a Theorem of Veech \cite{veech}, there exists a probability measure 
$\Prob$ on $\omab$ such that
the dynamical systems $(\omab, \sigma, \Prob)$ and 
$(\Delta(\R), \gz, \nu)$ are isomorphic. 
Indeed, let $w_1\in\A$ and define the set $\Delta(w_1)$ in the following way.
If $w_1=(a, n_1, \pi_1)$, then
$$
\Delta(w_1)=\{(\la,\pi)\in \Delta^- | \
\exists \lap\in\Delta_{a^{n_1}\pi}^+ {\rm \ such \ that} \ \la=\frac{A(w_1)\lap}{|A(w_1)\lap|}\}.
$$
If $w_1=(b, n_1, \pi_1)$, then
$$
\Delta(w_1)=\{(\la,\pi)\in \Delta^+ | \
\exists \lap\in\Delta_{b^{n_1}\pi}^- {\rm \ such \ that} \ \la=\frac{A(w_1)\lap}{|A(w_1)\lap|}\}.
$$

In other words, for a letter $w_1=(c_1, n_1, \pi_1)$,
the set $\Delta(w_1)$ is the set of all interval exchanges $(\la,\pi)$
such that, first,  $\pi=\pi_1$, and, second, that 
the application of the map $\gz$ to $(\la,\pi)$
results in $n_1$ applications of the Rauzy operation $c_1$.

The coding map $\Phi:\Delta(\R)\to\Omega_{\A, B}$ is given by the formula
\begin{equation}
\Phi(\la,\pi)=\omega_1 \dots \omega_n \dots \ {\rm if} \ 
\gz^n(\la,\pi)\in\Delta(\omega_n).
\end{equation}

The $\gz$-invariant smooth probability measure $\nu$ projects under
$\Phi$ to a $\sigma$-invariant measure
on $\Omega_{\A, B}$; probability with respect to
that measure will be denoted by $\Prob$.

For $w\in W_{\A,B}$, $w=w_1\dots w_n$, let
$$
C(w)=\{\omega\in\Omega_{\A,B} | \ \omega_1=w_1, \dots, \omega_n=w_n\}.
$$
and
$$
\Delta(w)=\Phi^{-1}(C(w)).
$$

W. Veech \cite{veech} has proved the following
\begin{proposition}
Let $\q\in\wab$ be a word such that all entries of the
matrix $A(\q)$ are positive.
Let $\omega\in\omab$ be a sequence having infinitely many occurrences of the word $\q$.
Then the set $\Phi^{-1}(\omega)$ consists of one point.
\end{proposition}

We thus obtain an almost surely bijective symbolic dynamics for the map $\gz$.

\subsection{Veech's space of zippered rectangles.}

A {\it zippered rectangle} associated to the Rauzy class $\R$ is a quadruple $(\la, h, a, \pi)$,  where $\la \in 
{\mathbb R}_+^m, 
h\in{\mathbb R}^m_+, a\in{\mathbb R}^m, \pi\in{\cal R}$, and the vectors 
$h$ and $a$ satisfy the following equations and inequalities (one introduces 
auxiliary components $a_0=h_0=a_{m+1}=h_{m+1}=0$, and sets $\pi(0)=0$, 
$\pi^{-1}(m+1)=m+1$): 

\begin{equation}
\label{zipone}
h_i-a_i=h_{\pi^{-1}(\pi(i)+1)}-a_{\pi^{-1}(\pi(i)+1)-1}, i=0, \dots, m
\end{equation}
\begin{equation}
\label{ziptwo}
h_i\geq 0,   i=1, \dots,m, \  
a_i\geq 0,  i=1, \dots, m-1,  
\end{equation}
\begin{equation}
\label{zipthree}
a_i\leq \min (h_i, h_{i+1}) {\rm \ for\ } i\neq m,i\neq \pi^{-1}m,  
\end{equation}
\begin{equation}
\label{zipfour}
a_m\leq h_m, \ a_m\geq -h_{\pi^{-1}m},\ 
a_{\pi^{-1}m}\leq h_{\pi^{-1}m+1} 
\end{equation}

The area of a zippered rectangle is given by the expression 
$$
\la_1h_1+\dots+\la_mh_m.
$$
Following Veech, we denote by $\Omega(\R)$ the space of all zippered 
rectangles, corresponding to a given Rauzy class $\R$ and satisfying the condition 
$$
\la_1h_1+\dots+\la_mh_m=1.
$$
We shall denote 
by $x$ an individual zippered rectangle.

Veech \cite{veech} defines a flow $P^t$ and a map $\U$ on the space of zippered
rectangles by the formulas:
$$
P^t(\la,h,a,\pi)=(e^t\la, e^{-t}h, e^{-t}a, \pi).
$$
$$
\U(\la,h,a,\pi)=\begin{cases} (A^{-1}(a,\pi)\la, A^t(a,\pi)h, a^{\prime}, a\pi), &\text{if $(\la,\pi)\in\Delta^-$}\\
 (A^{-1}(b,\pi)\la, A^t(b,\pi)h, a^{\prime\prime}, b\pi),&\text{if $(\la,\pi)\in\Delta^+$},
\end{cases}
$$
where
$$
a^{\prime}_i=\begin{cases} a_i, &\text{if $j<\pi^{-1}m$,}\\
h_{\pi^{-1}m}+a_{m-1}, &\text {if $i=\pi^{-1}m$,}\\
a_{i-1}, &\text{other $i$}.
\end{cases}
$$
$$
a^{\prime\prime}_i=
\begin{cases} a_i, &\text{if $j<m$,}\\
-h_{\pi^{-1}m}+a_{\pi^{-1}m-1}, &\text {if $i=m$.}
\end{cases}
$$

The map $\U$ is invertible; $\U$ and $P^t$ commute (\cite{veech}).
Denote
$$
\tau(\lambda, \pi)=(\log(|\la|-\min(\la_m,\la_{\pi^{-1}m})),
$$
and for $x\in \Omega(\R)$, $x=(\la, h, a ,\pi)$, write
$$
\tau(x)=\tau(\la,\pi).
$$
Now, following Veech \cite{veech}, define
$$
\Y(\R)=\{x\in\Omega(\R)\; |\; |\la|=1\}.
$$
and
$$
\Omega_0(\R)=\bigcup_{x\in\Y(\R), 0\leq t\leq \tau(x)}P^tx.
$$
$\Omega_0(\R)$ is a fundamental domain for $\U$ and,
identifying the points $x$ and $\U x$ in $\Omega_0(\R)$, we obtain
a natural flow, also denoted by $P^t$, on $\Omega_0(\R)$.
The space $\Omega_0(\R)$ will be referred to as {\it Veech's moduli space of zippered rectangles}, 
and the flow $P^t$ as the Teichm{\"u}ller flow on the space of zippered rectangles.

The space $\Omega(\R)$ has a natural Lebesgue measure class and so does
the transversal $\Y(\R)$.
Veech \cite{veech} has proved the following Theorem.

\begin{theorem}[Veech \cite{veech}]
\label{zipmeas}
There exists a measure $\mu_{\R}$ on $\Omega(\R)$, absolutely continuous with respect to Lebesgue, preserved by both the map
$\U$  and the flow $P^t$ and such that $\mu_{\R}(\Omega_0(\R))<\infty$.
\end{theorem}

\subsection{Symbolic representation for the flow $P^t$}

Following Zorich \cite{zorich},  define
$$
\Omega^+(\R)=\{ x=(\la, h,a,\pi)| \ (\la,\pi)\in\Delta^+, a_m\geq 0 \}.
$$
$$
\Omega^-(\R)=\{ x=(\la, h,a,\pi)| \ (\la,\pi)\in\Delta^-, a_m\leq 0 \},
$$
$$
\Y^+(\R)=\Y(\R)\cap \Omega^+(\R),\  \Y^-(\R)=\Y(\R)\cap \Omega^-(\R),\  \Y^{\pm}(\R)=\Y^+(\R)\cup \Y^-(\R).
$$
Take $x\in\Y^{\pm}(\R)$,  $x=(\la,h,a,\pi)$, and let
$\F(x)$ be the first return map of the flow $P^t$ on the transversal $\Y^{\pm}(\R)$.
The map $\F$ is a lift of the map $\gz$ to the space of zippered rectangles: 
\begin{equation}
{\rm if }\ \F(\la,h,a,\pi)=(\lap, h^{\prime}, a^{\prime}, \pip), \ {\rm then} \  (\lap, \pip)=\gz(\lap, \pip).
\end{equation}

Note that if $x\in\Y^+$, 
then $\F(x)\in \Y^-$, and if $x\in\Y^-$, then $\F(x)\in
\Y^+$).

The map $\F$ preserves a natural absolutely continuous invariant measure on
$\Y^{\pm}(\R)$: indeed,
since $\Y^{\pm}(\R)$ is a transversal to the flow $P^t$,
the measure $\mu_{\R}$ induces an absolutely
continuous measure ${\overline \nu}$ on $\Y^{\pm}(\R)$;
since $\mu_{\R}$ is both $\U$
and $P^t$-invariant, the measure ${\overline \nu}$ is $\F$-invariant.
Zorich \cite{zorich} proved that the measure ${\overline \nu}$ is finite
and ergodic for $\F$.

There exists a probability measure $\Prob$ on $\omabz$ such that the
 dynamical system $(\Y^{\pm}, {\overline \nu}, \F)$ is
measurably isomorphic to the system $(\omabz, \Prob, \sigma)$
\cite{veech, bufetov}.
The Teichm{\"u}ller flow $P^t$ on the space
$\Omega_0(\R)$ is a suspension flow over the map $\F$ on $\Y^{\pm}$.
Identifying $\Y^{\pm}$ and $\omabz$, we obtain a symbolic dynamics
for the Teichm{\"u}ller flow in the space of zippered rectangles.

\subsection{Formulation of the main result.}

We now fix a Rauzy class $\R$ and consider the space $\Delta(\R)$ 
of interval exchange transformations corresponding to $\R$.  
We denote by $\gz$ the Rauzy-Veech-Zorich induction map on $\R$
(for detailed definition, see \cite{bufetov}). 
As above, $\nu$ stands for the Zorich invariant measure 
of Theorem \ref{zorichthm}. 
The map $\gz$ is not mixing: there exists a decomposition  
$$
\Delta(\R)=\Delta^+(\R)\bigcup \Delta^-(\R),
$$
such that $\gz(\Delta^+(\R))=\Delta^-(\R)$ and  $\gz(\Delta^-(\R))=\Delta^+(\R)$.

The restriction of the map $\gz^2$ on $\Delta^+(\R)$ is, however, 
exact \cite{bufetov}, and, consequently, strongly mixing. 

We  endow the space $\Delta(\R)$ with the Hilbert metric $d$
by lifting the standard Hilbert metric on a simplex to $\Delta(\R)$ 
(for a formal definition, see  \cite{bufetov}), and, 
for an $\alpha$-H\"older function $\phi:\Delta(\R)\to{\mathbb R}$, by
$|\phi|_{H_{\alpha}}$ we denote its $\alpha$-H\"older norm.
 
The main result of this note is the following.

\begin{theorem}
\label{main}
Let $\gz:\Delta(\R)\to\Delta(\R)$ be the Rauzy-Veech-Zorich 
induction map and let $\nu$ be the absolutely continuous  
invariant measure.  Let $p>2$.
Then, for any $\alpha>0$, there exist positive constants $C, \delta$ such that 
for any  $\phi\in H_{\alpha}\cap L_p(\Delta^+(\R), \nu)$ and $\psi\in L_2(\Delta^+(\R), \nu)$ we have 
$$
|\int \phi\ \times \psi\circ \gz^{2n} d\nu-\int\phi d\nu\int\psi d\nu|\leq C\exp(-\delta n)
(|\phi|_{H_{\alpha}}+|\phi|_{L_p})(|\psi|_{L_2}).  
$$
\end{theorem}

{\bf Remark.} A stretched-exponential estimate has been obtained in \cite{bufetov}.

\section{Proof of Theorem \ref{main}.}

\subsection{The tower method.}
The proof proceeds by the tower method of Lai-Sang Young \cite{lsy}. 

First, following Veech, we 
consider a word $\q$ such that all entries of the renormalization matrix
$A(\q)$ are positive. 
An observation, going back to 
Veech \cite{veech}, states that the first return 
map of $\gz$ on $\Delta_{\q}$ of $\Delta(\R)$ is uniformly expanding.

Now, for any $N\in {\mathbb N}$, introduce the set 
$$
B(N,\q)=\{(\la,\pi)\in \Delta_{\q}| 
\gz^n(\la, \pi)\notin\Delta_{\q} \ {\rm \ for \ all}\ n=1, \dots, N\}.
$$

\begin{lemma}
\label{expreturn}
There exist positive constants $C, \theta$, depending only on $\q$ and such 
that for any $N$ we have 
$$
\Prob(B(N,\q))\leq C\theta^N.
$$
\end{lemma}

In view of the Lai-Sang Young Theorem \cite{lsy}, Lemma \ref{expreturn} implies 
Theorem \ref{main} (note that the ``uniform'' version presented here  
follows from the statement in \cite{lsy} by the Banach-Steinhaus Theorem).

It remains to establish Lemma \ref{expreturn}.
For the Teichm{\"u}ller flow, 
exponential estimates for return times of the flow into compact sets already 
exist (\cite{athreya}, \cite{bufetov}, \cite{agy}).
In what follows, exponential estimate for the map are derived from those of the flow.  
More precisely, we shall  compare the {\it continuous} return times for the flow
with the {\it discrete} times for the induction map.

\subsection{Reduction to estimates on return times for the Teichm{\"u}ller flow.}

Generally speaking, the continuous and the discrete return times are not comparable  
if the suspension function of the flow is bounded neither from above nor from below. 
Nonetheless, for return times into a set of the form $\Delta_{\q}$, where $q$ is a word such that 
all entries of the corresponding renormalization matrix are positive, we shall show that, indeed, 
the continuous and the discrete return times are comparable up to a multiplicative and an additive 
constant. 
This comparison will follow from Corollary 9 in \cite{bufetov}, which 
states that  the norm of the Rauzy-Veech-Zorich 
renormalization matrix $A(w)$ grows exponentially in the length of 
a word $w\in\wab$ and will allow us to conclude the proof of the Theorem.  
We proceed to precise statements.

\begin{lemma}
\label{comparison}
There exists a positive 
constant $C(\R)$, depending only on the Rauzy class $\R$ such that
for any $w\in\wab$ we have  
$$
\frac{|w|}{\log ||A(w)||}\leq C(\R).
$$
\end{lemma}

This statement immediately follows from Corollary 9 in \cite{bufetov} 
[it suffices to take the logarithm].
Now take $(\la,\pi)\in\Delta(\R)$
and consider the corresponding symbolic 
expansion 
$$
\Phi(\la,\pi)=(\omega_1, \dots, \omega_n, \dots).
$$
Set  
$$
n_{\q}(\la,\pi)=\min \{n>0| \gz^n(\la,\pi)\in\Delta_{\q}).
$$

Introduce the word $w_{\q}(\la,\pi)$ 
by the formula  
$$
w=\omega_1\dots \omega_{n_{\q}}
$$
and denote 
$$
\eta_{\q}(\la,\pi)=\log ||A(w_{\q}(\la,\pi))||.
$$
and 
$$
\tau_{\q}(\la,\pi)=\log |A(w_{\q}(\la,\pi))\la|.
$$

By definition, $\tau_q(\la,\pi)$ is the time it 
takes a zippered rectangle with base $(\la,\pi)$ 
under the Teichm{\"u}ller flow to reach $\Delta_q$.
Note that, by definition,  $\eta_q>\tau_q$. 
Generally speaking, these quantities are not comparable: 
for some interval exchanges $(\la,\pi)$, it may happen that 
$\tau_{\q}(\la,\pi)$ is very small while $\eta_{\q}(\la,\pi)$ may be large. 
Nonetheless, we have

\begin{proposition}
\label{morecomparison}
There exists a constant $C_{11}(q)$ such that  for all  $(\la, \pi)\in\Delta_{q}$
we have 
$$
\eta_{\q}(\la,\pi)-\tau_{\q}(\la,\pi)<C_{11}(q)
$$ 
\end{proposition}

This follows from the fact that 
if $(\la, \pi)\in\Delta_{\q}$, then there exists a 
constant $C_{21}(\q)$ such that 
$$
\frac{\la_i}{\la_j}\leq C_{21}(\q),
$$

and then for {\it any} matrix $A$ with nonnegative entries we have 

$$
\frac{|A\la|}{||A||}\geq C_{31}(\q)
$$

for some constant $C_{31}(\q)$, depending only on $\q$, and the 
Proposition is proved.

We now complete the Proof of Lemma \ref{expreturn}.
There exists $\epsilon>0$, depending only on $\q$ and such 
that we have  the following estimate
\begin{equation}
\label{intest}
\int_{\Delta} \exp(\epsilon \tau_{\q})d\Prob<\infty.
\end{equation}

This estimate, extending earlier estimates of J. Athreya on returns of 
the Teichm{\"u}ller flow into compact sets,  
is obtained in \cite{bufetov}, \S 11 
and also, independently, in \cite{agy}, \S 6. 

Proposition \ref{morecomparison} now yields
\begin{equation}
\label{intesttwo}
\int_{\Delta} \exp(\epsilon \eta_{\q})d\Prob<\infty, 
\end{equation}

which, by Lemma \ref{comparison}, implies, for some ${\tilde \epsilon}$, 
depending only on $\q$, the estimate:

\begin{equation}
\label{intestthree}
\int_{\Delta} \exp({\tilde \epsilon} n_{\q})d\Prob<\infty, 
\end{equation}
whence, finally, we obtain
$$
\Prob(B(N,\q))\leq C\theta^N,
$$
the result of Lemma \ref{expreturn}.
Theorem \ref{main} is proved completely.

{\noindent \bf Acknowledgements.} 
We are deeply grateful to Giovanni Forni for our 
many stimulating discussions and to the referee for the careful 
reading of the manuscript and suggestions on improving the presentation.
We gratefully acknowledges the hospitality of the Institut 
Henri Poincar{\'e} in Paris and Institut des Math{\'e}matiques 
in Luminy, where part of this work was done. 
In summer 2005 A.B. was a Clay Liftoff Fellow. 
A.B.'s research is partially supported by 
the National Science Foundation under grant DMS0604386.  This research was
partially conducted during the period A.A. served as a Clay Research
Fellow.

\end{document}